\documentclass[twoside,reqno]{amsart}
\usepackage{latexsym}%%triangles:\lhd\rhd\unlhd\unrhd%
\newcommand{\qdn}{\hspace*{-1.5mm}}
\newcommand{\qqdn}{\hspace*{-2.5mm}}
\newcommand{\xqdn}{\hspace*{-5.0mm}}
\newcommand{\xxqdn}{\hspace*{-10mm}}

\newcommand{\binm}{\binom}

\newcommand{\nnm}{\nonumber}
\newcommand{\be}{\begin{equation}}
\newcommand{\ee}{\end{equation}}
\newcommand{\ba}{\begin{array}}
\newcommand{\ea}{\end{array}}
\newcommand{\bmn}{\begin{eqnarray}}
\newcommand{\emn}{\end{eqnarray}}
\newcommand{\bnm}{\begin{eqnarray*}}
\newcommand{\enm}{\end{eqnarray*}}
\newcommand{\bln}{\begin{subequations}}
\newcommand{\eln}{\end{subequations}}

\newtheorem{thm}{Theorem}%[section]

\newtheorem{entry}{Entry}%%%%%%%%%%%%%%%%

\newcommand{\bbtm}[4]{\bibitem{kn:#1}{#2,}~{#3,}~{#4.}}
\newcommand{\cito}[1]{\cite{kn:#1}}
\newcommand{\citu}[2]{\cite[#2]{kn:#1}}

\begin{document} %%%%%%%%%% This paper is published in %%%%%%%
{%%\fns \today\hfill\copyright%% Printed in China} %%%%%%%%%%%%%%%
%%%%%%%%%%%%%%%%%%%%%%%%%%%%%%%%%%%%%%%%%%%%%%%%%%%%%%%%%%%%%%
\title{Chu-Vandermonde convolution and \\harmonic number identities}
\author{$^A$Chuanan Wei, $^B$Dianxuan Gong, $^A$Qin Wang}

\footnote{\emph{2010 Mathematics Subject Classification}: Primary
05A19 and Secondary 33B15, 33C20}

\dedicatory{$^A$Department of Information Technology\\
 Hainan Medical College, Haikou 571101, China\\
 $^B$College of Sciences\\
 Hebei Polytechnic University, Tangshan 063009, China}

 %Email address: \emph{}}
\thanks{\emph{Email addresses}: weichuanan@yahoo.com.cn (C. Wei),
gongdianxuan@yahoo.com.cn (D. Gong)}

 \keywords{Chu-Vandermonde convolution;
 Derivative operator; Harmonic number identity}

\begin{abstract}
By applying the derivative operators to Chu-Vandermonde convolution,
several general harmonic number identities are established.
\end{abstract}

%%%%%%%%%%%%%%%%%%%%%%%%%%%%%%%%%%%%%%%%%%%%%%%%%%%%%%%%%%%%%%%%%%%
\maketitle\thispagestyle{empty}%%%%%%%%%%%%%%%%%%%%%%%%%%%%%%%%%%%%
\markboth{C. Wei,  D. Gong, Q. Wang}%%%%%%%%%%%%%%%%%%%%%%%%%%%%
         {Chu-Vandermonde convolution and harmonic number identities}

%%%%%%%%%%%%%%%%%%%%%%%%%%%%%%%%%%%%%%%%%%%%%%%%%%%%%%%%%%%%%%%%%%%
%%%%%%%%%%%%%%%%%%%%%%%%%%%%%%%%%%%%%%%%%%%%%%%%%%%%%%%%%%%%%%%%%%%
%%%%%%%%%%%%%%%%%%%%%%%%%%%%%%%%%%%%%%%%%%%%%%%%%%%%%%%%%%%%%%%%%%%
\section{Introduction}
%%%%%%%%%%%%%%%%%%%%%%%%%%%%%%%%%%%%%%%%%%%%%%%%%%%%%%%%%%%%%%%%%%%
%%%%%%%%%%%%%%%%%%%%%%%%%%%%%%%%%%%%%%%%%%%%%%%%%%%%%%%%%%%%%%%%%%%

For $\l\in \mathbb{N}$ and $n\in \mathbb{N}_0$, define the
generalized harmonic numbers by
 \[H_{0}^{\langle l\rangle}=0\quad\text{and}\quad
  H_{n}^{\langle l\rangle}
  =\sum_{k=1}^n\frac{1}{k^l}\quad\text{for}\quad n=1,2,\cdots.\]
When $l=1$, they reduce to the classical harmonic numbers:
\[H_{0}=0\quad \text{and}\quad H_{n}
=\sum_{k=1}^n\frac{1}{k}\quad \text{for}\quad n=1,2,\cdots.\]
 There exist many elegant identities involving harmonic numbers. They
can be found in the papers \cito{andrews}, \cito{chen},
\cito{chu-a}, \cito{chu-b}, \cito{liu}, \cito{paule}, \cito{sofo-a},
\cito{sofo-b} and \cito{zheng}.

For two differentiable functions $f(x)$ and $g(x,y)$, define
respectively the derivative operator $\mathcal{D}_x$ and
$\mathcal{D}^2_{xy}$ by
 \bnm
&&\mathcal{D}_xf(x)=\frac{d}{dx}f(x)\Big|_{x=0},\\[1mm]
&&\mathcal{D}_{xy}^2g(x,y)=\frac{\partial^2}{\partial x\partial
y}g(x,y)\Big|_{x=y=0}.
 \enm

 Then it is not difficult to show the following two
 derivatives:
 \bnm
&&\mathcal{D}_x\:\binm{s+x}{t}=\binm{s}{t}\big\{H_s-H_{s-t}\big\},\\
&&\mathcal{D}_{xy}^2\:\binm{s+x}{t}\binm{u+y}{v}
=\binm{s}{t}\binm{u}{v}\big\{H_s-H_{s-t}\big\}\big\{H_u-H_{u-v}\big\},
 \enm
where $s,t,u,v\in N_0$ with $t\leq s$ and $v\leq u$.

There are numerous binomial identities in the literature. Thereinto,
Chu-Vandermonde convolution (cf. \citu{andrews-r}{p. 67}) can be
stated as
 \bmn
 &&\xxqdn\xqdn\sum_{k=0}^n\binm{x}{k}\binm{y}{n-k}=\binm{x+y}{n}. \label{vandermonde}
 \emn
By applying the derivative operators $\mathcal{D}_x$ and
$\mathcal{D}_{xy}^2$ to \eqref{vandermonde}, several general
harmonic number identities will be established in the next section.

%%%%%%%%%%%%%%%%%%%%%%%%%%%%%%%%%%%%%%%%%%%%%%%%%%%%%%%%%%%%%%%%%%%
%%%%%%%%%%%%%%%%%%%%%%%%%%%%%%%%%%%%%%%%%%%%%%%%%%%%%%%%%%%%%%%%%%%
\section{Harmonic number identities}
%%%%%%%%%%%%%%%%%%%%%%%%%%%%%%%%%%%%%%%%%%%%%%%%%%%%%%%%%%%%%%%%%%%

%%%%%%%%%%%%%%%%%%%%%%%%%%%%%%%%%%%%%%%%%%%%%%%%%%%%%%%%%%%%%%%%%%%
\subsection{}
%%%%%%%%%%%%%%%%%%%%%%%%%%%%%%%%%%%%%%%%%%%%%%%%%%%%%%%%%%%%%%%%%%%

Performing the replacements $x\to -p-1-x$, $y\to -q-1-y$ for
\eqref{vandermonde} with $p,q\in N_0$, we obtain the binomial sum:
  \bmn \label{vandermonde-a}
\sum_{k=0}^n\binm{p+k+x}{k}\binm{q+n-k+y}{n-k}=\binm{p+q+n+1+x+y}{n}.
 \emn
Applying the derivative operator $\mathcal{D}_x$ to
\eqref{vandermonde-a} and then letting $y=0$, we have
\bmn\label{harmonic-a}
 \sum_{k=0}^n\binm{p+k}{p}\binm{q+n-k}{q}H_{p+k}
=\binm{p+q+n+1}{n}\Big\{H_{p}+H_{p+q+n+1}-H_{p+q+1}\Big\}.
 \emn
When $p=q=0$, \eqref{harmonic-a} reduces to the known result (cf.
\citu{chen}{Equation 2.1}):
 \bnm
 \sum_{k=0}^nH_k=(n+1)\Big\{H_{n+1}-1\Big\}.
 \enm
Employing the substitutions $p\to p+1$, $k\to k-1$, $n\to n-1$ for
 \eqref{harmonic-a}, we have
 \bmn
\sum_{k=1}^n\binm{p+k}{k}\binm{q+n-k}{q}
 kH_{p+k}&&\xqdn=\:(p+1)\binm{p+q+n+1}{n-1} \nnm\\\label{harmonic-b}
&&\xqdn\times\:\Big\{H_{p+1}+H_{p+q+n+1}-H_{p+q+2}\Big\}.
 \emn
When $p=q=0$, \eqref{harmonic-b} reduces to the known result (cf.
\citu{chen}{Equation 2.2}):
 \bnm
 \sum_{k=1}^nkH_k=\frac{n(n+1)}{2}H_n-\frac{(n-1)n}{4}.
 \enm
The method, which is used to deduce \eqref{harmonic-b}, can further
be utilized to derive equations with the factor $k^m$, where
$m\geq2$. Considering that the general resulting identities will
become complicated, we shall only recover the two known results (cf.
\citu{chen}{Equations 2.3-2.4}):
 \bnm
&&\xqdn\sum_{k=1}^nk^2H_k=\frac{n(n+1)(2n+1)}{6}H_n-\frac{(n-1)n(4n+1)}{36},\\
&&\xqdn\sum_{k=1}^nk^3H_k=\frac{n^2(n+1)^2}{4}H_n-\frac{(n-1)n(n+1)(3n-2)}{48}.
 \enm

Applying the derivative operator $\mathcal{D}_{xy}^2$ to
\eqref{vandermonde-a} and using \eqref{harmonic-a}, we establish the
 theorem.

\begin{thm}\label{thm-a}
For $p,q\in N_0$, there holds the harmonic number identity:
 \bnm
&&\xxqdn\sum_{k=0}^n\binm{p+k}{p}\binm{q+n-k}{q}
 H_{p+k}H_{q+n-k}=\binm{p+q+n+1}{n}
\Big\{\big(H_{p+q+1}^{\langle 2\rangle}-H_{p+q+n+1}^{\langle
2\rangle}\big)
\\&&\qquad\qquad\qquad\:+\:\big(H_p-H_{p+q+1}+H_{p+q+n+1}\big)
\big(H_q-H_{p+q+1}+H_{p+q+n+1}\big)\Big\}.
 \enm
\end{thm}

When $p=q=0$, Theorem \ref{thm-a} reduces to the identity:
 \bnm
\sum_{k=0}^nH_kH_{n-k}=
 (n+1)\Big\{\big(1-H_{n+1}^{\langle 2\rangle}\big)+\big(H_{n+1}-1\big)^2\Big\}.
 \enm

Employing the substitutions $p\to p+1$, $k\to k-1$, $n\to n-1$ for
Theorem \ref{thm-a}, we found the theorem.

\begin{thm}\label{thm-b}
For $p,q\in N_0$, there holds the harmonic number identity:
 \bnm
&&\xqdn\qqdn\sum_{k=1}^n\binm{p+k}{p}\binm{q+n-k}{q}
 kH_{p+k}H_{q+n-k}=\binm{p+q+n+1}{n-1}
\Big\{\big(H_{p+q+2}^{\langle 2\rangle}-\:H_{p+q+n+1}^{\langle
2\rangle}\big)\\&&\qquad\qquad+\big(H_{p+1}-H_{p+q+2}+H_{p+q+n+1}\big)
\big(H_q-H_{p+q+2}+H_{p+q+n+1}\big)\Big\}(p+1).
 \enm
\end{thm}

When $p=q=0$, Theorem \ref{thm-b} reduces to the identity:
 \bnm
\sum_{k=1}^nkH_kH_{n-k}=\frac{n(n+1)}{2}
 \Big\{H_{n+1}^2-H_{n+1}^{\langle 2\rangle}-2H_{n+1}+2\Big\}.
 \enm
Further, we can deduce the following two identities:
 \bnm
&&\xqdn\qdn\sum_{k=1}^nk^2H_kH_{n-k} =\frac{n(n+1)(2n+1)}{6}
 \Big\{H_{n+1}^2-H_{n+1}^{\langle 2\rangle}
-\frac{13n+5}{3(2n+1)}H_{n+1}+\frac{71n+37}{18(2n+1)}\Big\},\\
&&\xqdn\qdn\sum_{k=1}^nk^3H_kH_{n-k} =\frac{n^2(n+1)^2}{4}
 \Big\{H_{n+1}^2-H_{n+1}^{\langle 2\rangle}
-\frac{7n+5}{3(n+1)}H_{n+1}+\frac{35n+37}{18(n+1)}\Big\}.
 \enm
%%%%%%%%%%%%%%%%%%%%%%%%%%%%%%%%%%%%%%%%%%%%%%%%%%%%%%%%%%%%%%%%%%%
\subsection{}
%%%%%%%%%%%%%%%%%%%%%%%%%%%%%%%%%%%%%%%%%%%%%%%%%%%%%%%%%%%%%%%%%%%
Performing the replacements $x\to p+n+x$, $y\to q+n+y$ for
\eqref{vandermonde} with $p,q\in N_0$, we get the binomial sum:
 \bmn \label{vandermonde-b}
\qqdn\qdn\sum_{k=0}^n\binm{p+n+x}{k}\binm{q+n+y}{n-k}=\binm{p+q+2n+x+y}{n}.
 \emn
Applying the derivative operator $\mathcal{D}_y$ to
\eqref{vandermonde-b} and then letting $x=0$, we have
 \bmn \label{harmonic-c}
\sum_{k=0}^n\binm{p+n}{k}\binm{q+n}{n-k}
 H_{q+k}
=\binm{p+q+2n}{n}\Big\{H_{q+n}+H_{p+q+n}-H_{p+q+2n}\Big\},
 \emn
which is a special case of \citu{chu-a}{Theorem 1.5}.

Applying the derivative operator $\mathcal{D}_{xy}^2$ to
\eqref{vandermonde-b} and using \eqref{harmonic-c}, we establish the
theorem.

\begin{thm}\label{thm-c}
For $p,q\in N_0$, there holds the harmonic number identity:
 \bnm
&&\xxqdn\sum_{k=0}^n\binm{p+n}{k}\binm{q+n}{n-k}
 H_{p+n-k}H_{q+k}=\binm{p+q+2n}{n}
\Big\{\big(H_{p+q+n}^{\langle 2\rangle}-H_{p+q+2n}^{\langle
2\rangle}\big)
\\&&\qquad+\:\big(H_{p+n}+H_{p+q+n}-H_{p+q+2n}\big)
\big(H_{q+n}+H_{p+q+n}-H_{p+q+2n}\big)\Big\}.
 \enm
\end{thm}

When $p=q=0$, Theorem \ref{thm-c} reduces to the known result due to
Chen and Chu \citu{chu-b}{Example 3}:
  \bnm
 \sum_{k=0}^n\binm{n}{k}^2H_kH_{n-k} =\binm{2n}{n}
 \Big\{\big(H_n^{\langle2\rangle}-H_{2n}^{\langle 2\rangle}\big)
 +\big(2H_n-H_{2n}\big)^2\Big\}.
 \enm

Employing the substitutions $q\to q+1$, $k\to k-1$, $n\to n-1$ for
Theorem \ref{thm-c}, we found the theorem.

\begin{thm}\label{thm-d}
For $p,q\in N_0$, there holds the harmonic number identity:
 \bnm
&&\xqdn\qdn\sum_{k=1}^n\binm{p+n}{k}\binm{q+n}{n-k}
 kH_{p+n-k}H_{q+k}=\binm{p+q+2n-1}{n-1}
\Big\{\big(H_{p+q+n}^{\langle 2\rangle}-H_{p+q+2n-1}^{\langle
2\rangle}\big)\\&&\quad+\big(H_{p+n-1}+H_{p+q+n}\!-\!H_{p+q+2n-1}\big)
\big(H_{q+n}+H_{p+q+n}\!-\!H_{p+q+2n-1}\big)\Big\}(p+n).
 \enm
\end{thm}

When $p=q=0$, Theorem \ref{thm-d} reduces to the identity:
  \bnm
 \sum_{k=1}^n\binm{n}{k}^2kH_kH_{n-k}&&\xqdn=\:\frac{n}{2}\binm{2n}{n}
 \Big\{\big(H_n^{\langle2\rangle}-H_{2n-1}^{\langle 2\rangle}\big)
 \\&&\xqdn+\:\big(2H_n-H_{2n-1}\big)\big(2H_n-H_{2n-1}-1/n\big)\Big\}.
 \enm

Further, we can derive the following two identities:
 \bnm
 &&\xxqdn\sum_{k=1}^n\binm{n}{k}^2k^2H_kH_{n-k}=\frac{n^3}{4n-2}\binm{2n}{n}
 \bigg\{\Big(H_n^{\langle2\rangle}-H_{2n-1}^{\langle 2\rangle}\Big)
 \\&&\xxqdn\:\:+\:\Big(2H_n-H_{2n-1}\Big)\Big(2H_n-H_{2n-1}
 -\frac{2n^2-1}{2n^3-n^2}\Big)-\frac{(n-1)(2n^2-2n+1)}{n^3(2n-1)^2}\bigg\},\\
&&\xxqdn\sum_{k=1}^n\binm{n}{k}^2k^3H_kH_{n-k}=\frac{n^3(n+1)}{8n-4}\binm{2n}{n}
 \bigg\{\Big(H_n^{\langle2\rangle}-H_{2n-1}^{\langle 2\rangle}\Big)
 \\&&\xxqdn\:\:+\:\Big(2H_n-H_{2n-1}\Big)\Big(2H_n-H_{2n-1}
 -\frac{2n^2+4n-4}{2n^3+n^2-n}\Big)-\frac{3(n-1)(2n^2-2n+1)}{n^2(n+1)(2n-1)^2}\bigg\}.
 \enm

%%%%%%%%%%%%%%%%%%%%%%%%%%%%%%%%%%%%%%%%%%%%%%%%%%%%%%%%%%%%%%%%%%%
\subsection{}
%%%%%%%%%%%%%%%%%%%%%%%%%%%%%%%%%%%%%%%%%%%%%%%%%%%%%%%%%%%%%%%%%%%

Performing the replacements $x\to-p-1-x$, $y\to q+n+y$ for
\eqref{vandermonde} with $p,q\in N_{0}$, we achieve the binomial
sum:
 \bmn
\sum_{k=0}^n(-1)^k\binm{p+k+x}{k}\binm{q+n+y}{n-k}
=(-1)^n\binm{p-q+x-y}{n}.\label{vandermonde-c}
 \emn
Applying the derivative operator $\mathcal{D}_x$ to
\eqref{vandermonde-c} and then letting $y=0$, we have
 \bmn \label{harmonic-d}
\quad\:\:\sum_{k=0}^n(-1)^k\binm{p+k}{k}\binm{q+n}{n-k}H_{p+k} =
\begin{cases}
 A_n,\:\text{for}\:\:p-q\geq n,\\
 B_n,\:\text{for}\:\:0\leq p-q< n,\\
 C_n,\:\text{for}\:\:p-q<0,
 \end{cases}
 \emn
where
 \bnm
&&\xxqdn A_n= (-1)^n\binm{p-q}{n}\big\{H_{p}+H_{p-q}-H_{p-q-n}\big\},\\
 &&\xxqdn B_n=
(-1)^{1+p-q}\frac{(p-q)!(q-p+n-1)!}{n!},\\
&&\xxqdn
C_n=(-1)^n\binm{p-q}{n}\big\{H_{p}+H_{q-p-1}-H_{q-p+n-1}\big\} .
 \enm
Applying the derivative operator $\mathcal{D}_y$ to
\eqref{vandermonde-c} and then letting $x=0$, we have
 \bmn \label{harmonic-e}
\quad\:\:\sum_{k=0}^n(-1)^k\binm{p+k}{k}\binm{q+n}{n-k}H_{q+k} =
\begin{cases}
 D_n,\:\text{for}\:\:p-q\geq n,\\
 E_n,\:\text{for}\:\:0\leq p-q< n,\\
 F_n,\:\text{for}\:\:p-q<0,
 \end{cases}
 \emn
where
 \bnm
&&\xxqdn D_n=(-1)^n\binm{p-q}{n}\big\{H_{q+n}+H_{p-q}-H_{p-q-n}\big\},\\
 &&\xxqdn E_n=(-1)^{1+p-q}\frac{(p-q)!(q-p+n-1)!}{n!},\\
&&\xxqdn
F_n=(-1)^n\binm{p-q}{n}\big\{H_{q+n}+H_{q-p-1}-H_{q-p+n-1}\big\}.
 \enm
We point out that \eqref{harmonic-d} can be given by
\citu{chu-a}{Theorem 1.1} and \eqref{harmonic-e} can be offered by
 \citu{chu-a}{Theorem 1.5}.

Applying the derivative operator $\mathcal{D}^2_{xy}$ to
\eqref{vandermonde-c} and using
\eqref{harmonic-d}-\eqref{harmonic-e}, we establish the theorem.

\begin{thm}\label{thm-e}
For $p,q\in N_0$, there holds the harmonic number identity:
 \bnm
&&\xqdn\sum_{k=0}^n(-1)^k\binm{p+k}{k}\binm{q+n}{n-k}
 H_{p+k}H_{q+k}=
\begin{cases}
 U_n,\:\text{for}\:\:p-q\geq n,\\
 V_n,\:\text{for}\:\:0\leq p-q< n,\\
 W_n,\:\text{for}\:\:p-q<0,
 \end{cases}
 \enm
where
  \bnm
U_n&&\xqdn=\:(-1)^n\binm{p-q}{n}\Big\{\big(H_{p-q-n}^{\langle
2\rangle}-H_{p-q}^{\langle2\rangle}\big)
 \\&&\xqdn+\:\big(H_p+H_{p-q}-H_{p-q-n}\big) \big(H_{q+n}+H_{p-q}-H_{p-q-n}\big)\Big\},\\
V_n&&\xqdn=\:(-1)^{1+p-q}\frac{(p-q)!(q-p+n-1)!}{n!}
\\&&\xqdn\times\:\Big\{H_p+H_{q+n}+2H_{p-q}-2H_{q-p+n}
+\frac{2}{q-p+n}\Big\},\\
W_n&&\xqdn=\:(-1)^n\binm{p-q}{n} \Big\{\big(H_{q-p-1}^{\langle
2\rangle}-\!H_{q-p+n-1}^{\langle 2\rangle}\big)
 \\&&\xqdn+\:\big(H_{p}+H_{q-p-1}-H_{q-p+n-1}\big) \big(H_{q+n}+H_{q-p-1}
-H_{q-p+n-1}\big)\Big\}.
 \enm
\end{thm}

When $p=q=0$ with $n>0$, Theorem \ref{thm-e} reduces to the
 identity:
  \bnm
 \sum_{k=0}^n(-1)^k\binm{n}{k}H_k^2=\frac{1}{n}
 \Big\{H_n-\frac{2}{n}\Big\}.
 \enm

Employing the substitutions $p\to p+1$, $q\to q+1$, $k\to k-1$,
$n\to n-1$ for Theorem \ref{thm-e}, we found the theorem.

\begin{thm}\label{thm-f}
For $p,q\in N_0$, there holds the harmonic number identity:
 \bnm
&&\xqdn\sum_{k=1}^n(-1)^k\binm{p+k}{k}\binm{q+n}{n-k}
 kH_{p+k}H_{q+k}=
\begin{cases}
 U_n^{*},\:\text{for}\:\:p-q\geq n-1,\\
 V_n^{*},\:\text{for}\:\:0\leq p-q<n-1,\\
 W_n^{*},\:\text{for}\:\:p-q<0,
 \end{cases}
 \enm
where
  \bnm
U_n^*&&\xqdn=\:(-1)^n(p+1)\binm{p-q}{n-1}\Big\{\big(H_{p-q-n+1}^{\langle
2\rangle}-H_{p-q}^{\langle2\rangle}\big)
 \\&&\xqdn+\:\big(H_{p+1}+H_{p-q}-H_{p-q-n+1}\big)\big(H_{q+n}+H_{p-q}-H_{p-q-n+1}\big)\Big\},\\
V_n^*&&\xqdn=\:(-1)^{p-q}(p+1)\frac{(p-q)!(q-p+n-2)!}{(n-1)!}
\\&&\xqdn\times\:\Big\{H_{p+1}+H_{q+n}+2H_{p-q}-2H_{q-p+n-1}
+\frac{2}{q-p+n-1}\Big\},\\
W_n^*&&\xqdn=\:(-1)^n(p+1)\binm{p-q}{n-1}
\Big\{\big(H_{q-p-1}^{\langle 2\rangle}-\!H_{q-p+n-2}^{\langle
2\rangle}\big)
 \\&&\xqdn+\:\big(H_{p+1}+H_{q-p-1}-H_{q-p+n-2}\big) \big(H_{q+n}+H_{q-p-1}
-H_{q-p+n-2}\big)\Big\}.
 \enm
\end{thm}

When $p=q=0$ with $n>1$, Theorem \ref{thm-f} reduces to the
identity:
  \bnm
 \sum_{k=1}^n(-1)^k\binm{n}{k}kH_k^2=\frac{1}{1-n}
 \Big\{H_n-\frac{n^2+3n-2}{n(n-1)}\Big\}.
 \enm

Further, we can deduce the following two identities:
 \bnm
 &&\xxqdn \sum_{k=1}^n(-1)^k\binm{n}{k}k^2H_k^2=\frac{n}{(n-1)(n-2)}
 \Big\{H_n-\frac{2n^3+n^2-11n+6}{n(n-1)(n-2)}\Big\},\\
 &&\xxqdn\sum_{k=1}^n(-1)^k\binm{n}{k}k^3H_k^2=\frac{(n+1)n}{(n-1)(n-2)(n-3)}
 \\&&\qquad\qquad\qquad\!\times\:\Big\{H_n-\frac{(3n^4-4n^3-32n^2+62n-15)n-6}{(n+1)n(n-1)(n-2)(n-3)}\Big\},
 \enm
where $n>2$ for the first equation and $n>3$ for the second
equation.

\textbf{Remark:} With the change of the parameters $p$ and $q$,
Theorems \ref{thm-a}-\ref{thm-f} can produce more interesting
harmonic number identities. We shall not lay them out one by one
because of the triviality of the work.
%%%%%%%%%%%%%%%%%%%%%%%%%%%%%%%%%%%%%%%%%%%%%%%%%%%%%%%%%%%%%%%%%%
%%%%%%%%%%%%%%%%%%%%%%%%%%%%%%%%%%%%%%%%%%%%%%%%%%%%%%%%%%%%%%%%%%%

%%%%%%%%%%%%%%%%%%%%%%%%%%%%%%%%%%%%%%%%%%%%%%%%%%%%%%%%%%%%%%%%%%%
%%%%%%%%%%%%%%%%%%%%%%%%%%%%%%%%%%%%%%%%%%%%%%%%%%%%%%%%%%%%%%%%%%%
%%%%%%%%%%%%%%%%%%%%%%%%%%%%%%%%%%%%%%%%%%%%%%%%%%%%%%%%%%%%%%%%%%%

\end{document}